\documentclass{ifacconf}

\usepackage{graphicx}      
\usepackage{natbib}        

\usepackage{mathtools,amsmath,amssymb}
\usepackage{xcolor}
\usepackage{xspace}
\usepackage{booktabs}
\usepackage{colortbl}
\makeatletter
\newcommand{\ccell}[3][]{%
  \kern-\fboxsep
  \if\relax\detokenize{#1}\relax
    \expandafter\@firstoftwo
  \else
    \expandafter\@secondoftwo
  \fi
  {\colorbox{#2}}%
  {\colorbox[#1]{#2}}%
  {#3}\kern-\fboxsep
}
\makeatother
\definecolor{cellgray}{gray}{0.9}
\usepackage{tikz}
\usepackage{pgfplots}

\usepackage{siunitx}
\newcommand{\R}{\ensuremath\mathbb{R}}

\newcommand{\T}{\ensuremath\mathsf{T}}

\newcommand{\calP}{\mathcal{P}}

\newcommand{\lineWidth}{1.2pt}

\definecolor{mycolor1}{rgb}{0.00000,0.44700,0.74100}
\definecolor{mycolor2}{rgb}{0.85000,0.32500,0.09800}
\definecolor{mycolor3}{rgb}{0.92900,0.69400,0.12500}
\definecolor{mycolor4}{rgb}{0.46600,0.67400,0.18800}
\definecolor{mycolor5}{rgb}{0.49400,0.18400,0.55600}

\definecolor{nicegreen}{rgb}{0.3,0.7,0.4}

\newcommand{\timeInt}{\mathbb{T}}

\newcommand{\state}{x}
\newcommand{\stateDim}{n}

\newcommand{\reduce}[1]{\widehat{#1}}
\newcommand{\stateRed}{\reduce{\state}}
\newcommand{\stateDimRed}{r}

\newcommand{\param}{\mu}
\newcommand{\paramSet}{\calP}
\newcommand{\paramDim}{p}

\newcommand{\encoder}{\rho}
\newcommand{\decoder}{\varphi}

\definecolor{cbsblue}{RGB}{68,119,170}
\definecolor{cbscyan}{RGB}{102,204,238}
\definecolor{cbsgreen}{RGB}{34,136,51}
\definecolor{cbsyellow}{RGB}{204,187,68}
\definecolor{cbsred}{RGB}{238,102,119}
\definecolor{cbspurple}{RGB}{170,51,119}
\definecolor{cbsgrey}{RGB}{187,187,187}

\pgfplotscreateplotcyclelist{approx}{%
line width = \lineWidth, color=cbsblue\\
line width = \lineWidth, color=cbscyan\\
line width = \lineWidth, color=cbsgreen\\
line width = \lineWidth, color=cbsyellow\\
line width = \lineWidth, color=cbsred\\
line width = \lineWidth, color=cbspurple\\
line width = \lineWidth, color=cbsgrey\\
}

\pgfplotscreateplotcyclelist{approxEverySecondDashed}{%
line width = \lineWidth, color=cbsblue\\
line width = \lineWidth, color=cbscyan, dashed\\
line width = \lineWidth, color=cbsgreen\\
line width = \lineWidth, color=cbsyellow, dashed\\
line width = \lineWidth, color=cbsred\\
line width = \lineWidth, color=cbspurple, dashed\\
line width = \lineWidth, color=cbsgrey\\
}

\pgfplotscreateplotcyclelist{error}{%
line width = \lineWidth, color=cbscyan\\
line width = \lineWidth, color=cbsgreen\\
line width = \lineWidth, color=cbsyellow\\
line width = \lineWidth, color=cbsred\\
line width = \lineWidth, color=cbspurple\\
line width = \lineWidth, color=cbsgrey\\
}
\newcommand{\abbr}[1]{\textsf{#1}\xspace}
\newcommand{\FOM}{\abbr{FOM}}		
\newcommand{\MOR}{\abbr{MOR}}		
\newcommand{\ROM}{\abbr{ROM}}		

\newcommand{\POD}{\abbr{POD}}       
\newcommand{\IVP}{\abbr{IVP}}		
\newcommand{\NNA}{\abbr{NNA}}		
\newcommand{\LNA}{\abbr{LNA}}		
\newcommand{\NLA}{\abbr{NLA}}		
\newcommand{\LNAwt}{\textsf{LNA}$_{\withTime}$\xspace}

\newcommand{\withTime}{\text{ext}}
\newcommand{\stateTimeExt}{z}
\newcommand{\fTimeExt}{f_{\withTime}}
\newcommand{\VTimeExt}{V_{\withTime}}
\newcommand{\rhoTimeExt}{\rho_{\withTime}}
\newcommand{\decoderTimeExt}{\varphi_{\withTime}}
\newcommand{\WTimeExt}{W_{\withTime}}

\newcommand{\nrSnapshots}{M}
\newcommand{\LNAext}{$\text{\LNA}_{\text{ext}} $}
\newcommand{\LNAextTime}{$\text{\LNA}_{\text{ext,fix}} $}

\begin{document}
\begin{frontmatter}

\title{Leveraging time and parameters for nonlinear model reduction methods\thanksref{footnoteinfo}} 

\thanks[footnoteinfo]{BU is funded by Deutsche Forschungsgemeinschaft (DFG, German Research Foundation) under Germany's Excellence Strategy - EXC 2075 - 390740016 and by the BMBG (grant no.~05M22VSA). Further, BU acknowledges the support by the Stuttgart Center for Simulation Science (SimTech).}

\author[First]{Silke Glas} 
\author[Second]{Benjamin Unger} 

\address[First]{Department of Applied Mathematics, University of Twente, 7500 AE Enschede, The Netherlands (e-mail: s.m.glas@utwente.nl).}
\address[Second]{Stuttgart Center for Simulation Science (SC SimTech), University of Stuttgart, 70569 Stuttgart, Germany (e-mail: benjamin.unger@simtech.uni-stuttgart.de)}

\begin{abstract}                
In this paper, we consider model order reduction (MOR) methods for problems with slowly decaying Kolmogorov $n$-widths as, e.g., certain wave-like or transport-dominated problems. To overcome this Kolmogorov barrier within MOR, nonlinear projections are used, which are often realized numerically using autoencoders. These autoencoders generally consist of a nonlinear encoder and a nonlinear decoder and involve costly training of the hyperparameters to obtain a good approximation quality of the reduced system. To facilitate the training process, we show that extending the to-be-reduced system and its corresponding training data makes it possible to replace the nonlinear encoder with a linear encoder without sacrificing accuracy, thus roughly halving the number of hyperparameters to be trained. 
\end{abstract}


\end{frontmatter}

\section{Introduction}
In high-fidelity simulation-driven design processes, large-scale differential or differential-algebraic equations must be evaluated for many different sets of parameters, amounting to a significant demand on computing resources. Replacing the high-fidelity simulation with a cheap-to-evaluate surrogate model is expected to reduce the computing time and speed up the design cycle. Due to its efficient and numerically stable construction, generalization capabilities, and error certificates, \emph{model order reduction} (\MOR) is a typical approach for assembling the surrogate. Over the last decades, \MOR was successfully utilized in various application domains -- we refer to \cite{Ant05,BenCOW17} and the references therein for examples. 

In this paper, we consider projection-based \MOR methods, which construct the \ROM by projecting the differential equation onto a low-dimensional subspace using a Petrov--Galerkin framework. 
Although this approach is successful in many applications, transport-dominated problems with a strong coupling of the spatial and temporal domain typically require a rather large subspace to achieve good approximation quality. The theoretical limit is given by the Kolmogorov $n$-width \citep{Kol36}, respectively the Hankel singular values for control systems (see \cite{UngG19}). Different nonlinear approaches have been proposed in the literature to overcome this Kolmogorov barrier. Prominent examples are the shifted proper orthogonal decomposition \citep{ReiSSM18,BlaSU20,BlaSU22}, the registration method \citep{FerTZ22}, and methods based on (shallow) autoencoders \citep{LeeC20,KimCWZ22,BucGH21}. We refer to \cite{BucGHU24} for a more detailed literature review and a unifying framework for nonlinear projections (consisting of an encoder and a decoder when using autoencoders).

Most data-driven nonlinear approaches have in common that constructing the involved nonlinear projection entails solving complex and large-scale optimization problems. These are extremely expensive and require extensive fine-tuning to converge to a (good) local minimum. In addition, the theory-to-practice gap in deep learning \citep{GroV24} and the stability of the numerical method \citep{CohDPW22} prevent practical algorithms based on sampled data from achieving the theoretical optimum. 
Moreover, even if evaluated fast, these nonlinear approaches also require considerable computing power once they are deployed, see e.g., \cite{DesMH23}. 

Our main contribution is to simplify the structure of the nonlinear approach in the sense that the nonlinear projection can be computed with significantly fewer optimization parameters. Commonly, autoencoders consist of a nonlinear encoder and a nonlinear decoder. To achieve good accuracy, they involve a large number of hyperparameters. However, we demonstrate that extending the to-be-reduced system and its corresponding training data makes it possible to replace the nonlinear encoder with a linear encoder without sacrificing accuracy, thus roughly halving the number of hyperparameters to be trained.

The paper is structured as follows. We briefly introduce the setting and known work in Section~\ref{sec:knownWork}. Subsequently, we present the main idea of this paper, i.e., the extension of the training data, and some theoretical results in Section~\ref{sec:contribution}. In Section~\ref{sec:numExp}, we consider the numerical investigation of two proof-of-concept examples, the Burgers' equation and the advection equation. We conclude in Section~\ref{sec:conclusion}.

\section{Setting and known work}\label{sec:knownWork}

Consider the initial value problem (\IVP)
\begin{equation}
	\label{eqn:IVP}
	\left\{\quad \begin{aligned}
		\dot{\state}(t;\mu) &= f(t,\state(t;\param);\param),\\
		\state(t_0;\mu) &= \state_0(\param), 
	\end{aligned}\right.
\end{equation}
with state $\state(t;\param)\in\R^{\stateDim}$, set of parameters $\param\in\paramSet\subseteq\R^{\paramDim}$, and initial value $\state_0(\param)\in\R^{\stateDim}$. Our assumption is that $\stateDim$ is large and that the \IVP~\eqref{eqn:IVP} must be evaluated for many different $\param\in\paramSet$. The goal is thus to replace~\eqref{eqn:IVP} with a \emph{reduced-order model} (\ROM) of the form
\begin{equation}
	\label{eqn:ROM}
	\left\{\quad \begin{aligned}
		\dot{\stateRed}(t;\param) &= \reduce{f}(t,\stateRed(t;\param);\param),\\
		\stateRed(t_0;\param) &= \stateRed_0(\param), 
	\end{aligned}\right.
\end{equation}
with reduced state $\stateRed(t;\param)\in\R^{\stateDimRed}$, $\stateDimRed\ll\stateDim$ as an efficient surrogate of the parameter-to-solution map. In standard linear-subspace \MOR, the \ROM~\eqref{eqn:ROM} is constructed via a Petrov--Galerkin projection, i.e., one determines bi-orthogonal matrices $V,W\in\R^{\stateDim\times\stateDimRed}$ and constructs
\begin{align*}
	\reduce{f}(t,\stateRed;\param) \vcentcolon= W^\T f(t,V\stateRed;\param).
\end{align*}
When using linear-subspace \MOR methods, the best possible error that the \ROM can achieve is lower bounded by the Kolmogorov $\stateDim$-width.
To overcome this Kolmogorov barrier, the linear mappings encoded by the matrices $V$ and $W^\T$ are replaced by smooth nonlinear mappings
\begin{align*}
	\decoder&\colon\R^{\stateDimRed}\to\R^{\stateDim}, &
	\encoder&\colon\R^{\stateDim}\to\R^{\stateDimRed},
\end{align*}
which we refer to as the decoder and encoder, respectively. In the context of \MOR, these mappings have to satisfy the \emph{point projection property}, cf.~\cite{BucGHU24}, i.e.,
\begin{align}
	\label{eqn:projectionProperty}
	\encoder \circ \decoder = \mathrm{Id}_{\R^\stateDimRed},
\end{align}
where $\mathrm{Id}_{\R^\stateDimRed}$ denotes the identity mapping in $\R^{\stateDimRed}$. 
Allowing for this kind of nonlinear approximations, then, under certain smoothness assumptions, we have the following result taken from \cite[Cor.~3.6]{BucGHU24}, where also the exact setting and the proof can be found. 

\begin{thm}
	\label{thm:cor3-6}
	If the (time$\times$parameter)-to-solution map is sufficiently smooth, then there exists a \ROM of dimension $\paramDim+1$ that yields an error-free surrogate for~\eqref{eqn:IVP}.
\end{thm}

For practical purposes, the mappings $\encoder$ and $\decoder$ are obtained from snapshots, i.e., evaluations of~\eqref{eqn:IVP} for selected parameter and time values. In the following, the set of snapshots is denoted with $\{\state_i\}_{i=1}^{\nrSnapshots}$.

\section{Contribution}\label{sec:contribution}
The key observation for the proof of Theorem~\ref{thm:cor3-6} is, that the (time$\times$parameter)-to-solution map can serve as the decoder function, and the coordinates of the \ROM are simply the time-variable and the parameters. Nevertheless, in almost all linear-subspace and nonlinear \MOR methods, this information is not given during the training. Our key idea is thus to include the time (and the parameter) information in the training data. By doing so, we provide the autoencoder with additional information, which allows the decoder to (approximately) learn the (time$\times$parameter)-to-solution map from the given data. Subsequently, this means that the encoder can be realized by a linear mapping only, thus facilitating the training complexity. In more detail, we propose to perform an autonomization step (even if $f$ has no explicit time dependency), i.e., we consider the \emph{extended system} with respect to time and parameter
\begin{equation}
	\label{eqn:IVP:timeExtended}
	\left\{\quad \begin{aligned}
		\dot{\stateTimeExt}(t;\param) &= \fTimeExt(\stateTimeExt(t;\param);\param),\\
		\stateTimeExt(t_0;\param) &= \stateTimeExt_0(\param), 
	\end{aligned}\right.
\end{equation}
with extended state and vector field
\begin{align*}
	\stateTimeExt(t;\param) &\vcentcolon= \begin{bmatrix}
		t\\		
		\state(t;\param)\\
		\param
	\end{bmatrix}, &
	\fTimeExt(\stateTimeExt(t;\param);\param) &\vcentcolon= \begin{bmatrix}
		1\\
		f(t,\state(t;\param);\param)\\
		0
	\end{bmatrix}.
\end{align*}
The advantage of the formulation~\eqref{eqn:IVP:timeExtended} over~\eqref{eqn:IVP} is twofold:
\begin{enumerate}
	\item[1.] Using data from~\eqref{eqn:IVP:timeExtended} explicitly encodes the time variable and the parameters, which are now explicitly available in the \MOR process.
	\item[2.] There is no need for a nonlinear encoder function $\encoder$, since the time and the parameter can be extracted explicitly without the need for a nonlinear mapping.
\end{enumerate}

We denote the snapshots for~\eqref{eqn:IVP:timeExtended} corresponding to the snapshots $\{\state_i\}_{i=1}^{\nrSnapshots}$ with $\{\stateTimeExt_i\}_{i=1}^{\nrSnapshots}$.
As an immediate consequence of Theorem~\ref{thm:cor3-6}, we obtain the following result.

\begin{thm}
	\label{thm:linEncoder}
	If the (time$\times$parameter)-to-solution map is sufficiently smooth 
	 and $\stateDimRed = \paramDim + 1$, then
	\begin{align*}
	\begin{aligned}
	\min_{\substack{\encoder \colon  \R^{\stateDim}  \to  \R^{\stateDimRed} \\
	\decoder  \colon  \R^{\stateDimRed}  \to  \R^{\stateDim}  }} \sum_{i=1}^\nrSnapshots &\|\state_i - \decoder(\encoder(\state_i))\|_2^2\\
	 &= \min_{\substack{W\in \R^{(1+\stateDim+\paramDim)\times\stateDimRed}\\ \decoderTimeExt \colon \R^{\stateDimRed} \to  \R^{(1+\stateDim+\paramDim)} }} \smash{\sum_{i=1}^\nrSnapshots} \|\stateTimeExt_i - \decoderTimeExt(\WTimeExt^\T \stateTimeExt_i)\|_2^2,
	\end{aligned}
	\end{align*}
	with 
	\begin{align*}
	\WTimeExt &=
	\begin{bmatrix}
		W_{\mathrm{t}}, 
		W_{\mathrm{x}}, 
		W_{\mathrm{p}}
	\end{bmatrix} \in\R^{(1+\stateDim+\paramDim) \times  \stateDimRed}. 
\end{align*}	
\end{thm}

\begin{pf}
	We start by noticing, that the first term is equal to $0$ if we choose $\decoder$ as the (time$\times$parameter)-to-solution map and $\encoder= \decoder^{-1}$, due to Theorem~\ref{thm:cor3-6}. To show the equality to the extended system, we particularly choose $W_{\mathrm{t}} = 1$, $W_{\mathrm{x}} $ as the matrix with all entries 0 and $W_{\mathrm{p}}= [0_{(1+n)\times r}, \mathrm{Id}_{\R^p}]$, which leads to the coordinates of the \ROM being the time-variable and the parameters. Now setting 
	\begin{align*}
	\decoderTimeExt\colon  \R^{\stateDimRed} \to  \R^{(1+\stateDim+\paramDim)},\quad [t, \mu]^\T \mapsto [1, \decoder(t,\mu), 1_{1 \times p} ]^\T,
	\end{align*}
	 yields the desired result. 
\end{pf} 

\begin{rem}
	The use of an affine encoder as we propose it here is a common choice in \MOR based on quadratically or polynomially embedded manifolds, exemplified by \cite{Geelen2022,Barnett2022}.
Nevertheless, in these cases, the encoder is applied on the original state $\state$ and not on the extended state $\stateTimeExt$. Moreover, these approaches still suffer from (a modified) Kolmogorov barrier, which is due to the design choice of the decoder; see \cite{BucGH24} for details.
\end{rem}

Before we provide numerical results, we emphasize that it is solely possible to replace the encoder with a linear mapping, but not the decoder. If we happen to choose the decoder as a linear mapping we obtain the following result.

\begin{thm}
	\label{thm:linearDecoder}
	Consider snapshots $S = \{\state_i\}_{i=1}^{\nrSnapshots} \in\R^{\stateDim\times \nrSnapshots}$ and let $\sigma_i$ denote the $i$-th singular value of $S$. Then
	\begin{align}\label{eqn:linearDecoder}
		\min_{\substack{V\in\R^{\stateDim\times\stateDimRed}\\\encoder\colon\R^{\stateDim} \to \R^{\stateDimRed}}} \sum_{i=1}^M \|\state_i - V\encoder(\state_i)\|_2^2 = \sum_{j=r+1}^{\min(\nrSnapshots, \stateDim)}\sigma^2_{j}.
	\end{align}
Particularly, an optimal $V\in\R^{\stateDim\times\stateDimRed}$ is given by the leading~$\stateDimRed$ left singular vectors. Choosing a linear or nonlinear encoder doesn't affect the approximation error in this setting.
\end{thm}

\begin{pf}
	Without loss of generality, we can assume that the columns of $V$ are orthonormal. Then, the best approximation of each snapshot is given by the orthogonal projection onto the columns of $V$, i.e., for  given orthogonal $V\in\R^{\stateDim\times\stateDimRed}$,
	\begin{align*}
		\min_{\encoder\colon\R^{\stateDim}\to\R^{\stateDimRed}} \sum_{i=1}^M \|\state_i - V\encoder(\state_i)\|_2^2 = \sum_{i=1}^M \|\state_i - VV^\T \state_i\|_2^2.
	\end{align*}
	The claim follows from the Schmidt--Eckart--Young--Mirsky theorem, see, e.g., \cite[Thm.~3.6]{Ant05}.
\end{pf}

We emphasize that in contrast to Theorem~\ref{thm:cor3-6}, we do not distinguish between state $\state$ and extended state $\stateTimeExt$ in Theorem~\ref{thm:linearDecoder}. If we consider $\stateTimeExt$ in Theorem~\ref{thm:linearDecoder}, we get the best approximation error with respect to $\| \cdot \|_{\mathrm{F}}$ by
\begin{align*}
	\min_{\rhoTimeExt} \sum_{i=1}^M \|\stateTimeExt_i - \VTimeExt \rhoTimeExt(\stateTimeExt_i)\|_2^2 = \sum_{i=1}^M \|\stateTimeExt_i - \VTimeExt \VTimeExt^\T \stateTimeExt_i \|_2^2,
\end{align*}
with $\VTimeExt \in\R^{(1+\stateDim+\paramDim)\times\stateDimRed}$ and where we optimize over all $\rhoTimeExt\colon\R^{\stateDim+1}\to\R^{\stateDimRed}$.
Naturally, the resulting partial error in the state $\state$ in the latter equation can be (at best) as good as in \eqref{eqn:linearDecoder}, where we specifically chose the best approximation with respect to $\state$. This means, that we cannot expect any improvement in the state approximation accuracy by using the extended variable $\stateTimeExt$ only. Indeed, the following numerical example demonstrates this aspect. 

\begin{exmp}
	Consider $t_1 = 1$, $t_2 = 2$, and snapshots
	\begin{align*}
		\state_1 &= \begin{bmatrix}
			\alpha\\
			0
		\end{bmatrix}, & 
		\state_2 &= \begin{bmatrix}
			0\\
			\beta
		\end{bmatrix}, &
		\stateTimeExt_1 &= \begin{bmatrix}
			1\\
			\alpha\\
			0
		\end{bmatrix}, &
		\stateTimeExt_2 &= \begin{bmatrix}
			2\\
			0\\
			\beta
		\end{bmatrix},
	\end{align*}
	with $\alpha,\beta>0$. For $\alpha=1$ and $\beta = 0.9$, we obtain that for $\stateDimRed = 1$ the best approximation yields an error of $0.9$ in the Frobenius norm.
	Whereas, if we use the optimal approximation for the extended variable $z$, then the associated approximation of the state $x$ yields an error of approximately 0.98 with respect to $\| \cdot \|_{\mathrm{F}}$. This can also easily be seen via the interlacing property of the singular values \cite[Cor.~8.6.3]{gloub2013matrix}, as adding a row only enlarges the smallest singular value of a matrix. Moreover, if we split the matrix $\VTimeExt$ in the following way 
\begin{align*}
	\VTimeExt&=
	\begin{bmatrix}
		V_{\mathrm{t}}\\
		V_{\mathrm{x}}\\
		V_{\mathrm{p}}
	\end{bmatrix} 
	\in\R^{(1+\stateDim+\paramDim)\times\stateDimRed},
\end{align*}
and consider the projection error with respect to the state
\begin{align*}
	 \sum_{i=1}^{\nrSnapshots} \|\state_i - V_{\mathrm{x}} V_{\mathrm{x}}^\T(\state_i)\|_2^2,
\end{align*}
this value is of course also lower bounded by \eqref{eqn:linearDecoder}. However, this error also can be significantly larger, which is the case here, where the error measured in  $\| \cdot \|_{\mathrm{F}}$ is 1.26. 
\end{exmp}

\section{Numerical examples}\label{sec:numExp}
To demonstrate our theoretical findings, we consider two proof-of-concept examples, i.e., the Burgers' (Section~\ref{subsec:Burgers}) and the advection equation (Section~\ref{subsec:advection}). For both settings, we compare different autoencoder architectures, for which we describe the setup and the training process in Section \ref{subsec:autoenc_training}.

\subsection{Autoencoder Training} \label{subsec:autoenc_training}
We generate training data by sampling both examples, i.e., the Burgers' equation (Section~\ref{subsec:Burgers}) and the advection equation (Section~\ref{subsec:advection}). We pick a one-dimensional spatial domain $\Omega$ for both settings, which we discretize into $\stateDim=2^9 = 512$ equidistant points. For the time domain, we choose $300$ equidistant points and the generated simulation data is split into training, validation, and testing data by distributing the time points into $t_{3i}$, $t_{3i+1}$, and $t_{3i+2}$ for $i=1,\ldots,\nrSnapshots=100$ such that we have training data $\Theta\in\R^{\stateDim\times\nrSnapshots}$. We fix the parameter for both examples and only use time for the extended data. Consequently, the extended training data $\Theta_{\withTime}$ has dimension $(1+\stateDim)\times \nrSnapshots$.
With these data sets, we train five kinds of autoencoders: 
\begin{enumerate}
	\item[(I)] a \emph{nonlinear--nonlinear autoencoder} (\NNA) with nonlinear encoder and decoder trained on $\Theta$,
	\item[(II)] a \emph{linear-nonlinear autoencoder} (\LNAext) with linear encoder and nonlinear decoder trained on the time extended data $\Theta_{\withTime}$,
	\item[(III)] a \emph{linear-nonlinear autoencoder} (\LNA) with linear encoder and nonlinear decoder trained on~$\Theta$,
	\item[(IV)] a \emph{linear-nonlinear autoencoder} (\LNAextTime) with fixed linear encoder given by the first unit vector and nonlinear decoder trained on $\Theta_{\withTime}$, and
	\item[(V)] a \emph{nonlinear-linear autoencoder} (\NLA) with nonlinear encoder and linear decoder trained on data~$\Theta$. 
\end{enumerate}

For the implementation, we use \texttt{pytorch} \citep{NEURIPS2019_9015} and \texttt{Adam} \citep{kingma2014adam} with default settings as optimization algorithm. The target function is the least-squares error and  a batch size of 20 is set in all trainings. The optimization is terminated if there is no further improvement on the validation data after 100 epochs (with a maximum of 20000 epochs) to avoid overfitting. To reduce the impact of randomness in the model initialization and the training, we initialize and optimize each model 100 times and choose the  best run with respect to the lowest target function value on the validation error. 

To ensure comparability, we use the same architecture for all autoencoders: As an activation function, the Leaky ReLU is utilized. For the layers of the decoder, we consider the following three scenarios:
\begin{enumerate}
	\item[(A)] 7 layers with dimensions 1, 3, 9, 27, 81, 243, N,
	\item[(B)] 6 layers with dimensions 1, 4, 16, 64, 256, N,
	\item[(C)] 5 layers with dimensions 1, 5, 25, 125, N,
\end{enumerate}
where $N$ is either given by the spatial dimension, i.e., $N = \stateDim$ (for \NNA, \LNA, \NLA), or for the extended data by $N = 1+\stateDim$ (for \LNAext, \LNAextTime). If the encoder is nonlinear, we use the mirrored setting of the decoder. 

We report the minimum and the average error with respect to the testing data, with the error measure being the average over time of the relative $\| \cdot \|_2$-norm in space, i.e., 
\begin{align}\label{eqn:error_measure}
	\frac{1}{M} \sum_{i=1}^{M}  \frac{\| x_i - x_{i, \text{approx}}\|_2}{\| x_i \|_2}.
\end{align}

\subsection{Burgers' equation}
\label{subsec:Burgers}

\begin{figure}
	\begin{flushright}
	\begin{tikzpicture}
    \begin{axis}[
    width=.85\linewidth,
    height=1.2in,
    scale only axis,
    xmin=0, xmax=1,
    ymin=-0.05, ymax=0.65,
	grid=both,
	grid style={line width=.1pt, draw=gray!10},
	major grid style={line width=.2pt,draw=gray!50},
	axis lines*=left,
	axis line style={line width=\lineWidth},
    legend columns=4,
    legend entries={t=0.00, t=1.32, t=2.65, t=3.97},
    legend style={font=\footnotesize},
    legend cell align={left},
    legend pos = north west,
    cycle list name=approx,
    no markers,
    ]
    \addplot table[x=sampling, y=Ind-0,col sep=comma] {BurgersEquationLargeTime-FOM.dat};
    \addplot table[x=sampling, y=Ind-33,col sep=comma] {BurgersEquationLargeTime-FOM.dat};
    \addplot table[x=sampling, y=Ind-66,col sep=comma] {BurgersEquationLargeTime-FOM.dat};
    \addplot table[x=sampling, y=Ind-99,col sep=comma] {BurgersEquationLargeTime-FOM.dat};
\end{axis}
\end{tikzpicture}
	\end{flushright}
	\caption{Burgers' equation: solution at selected times.}
	\label{fig:BurgersData}
	\end{figure}
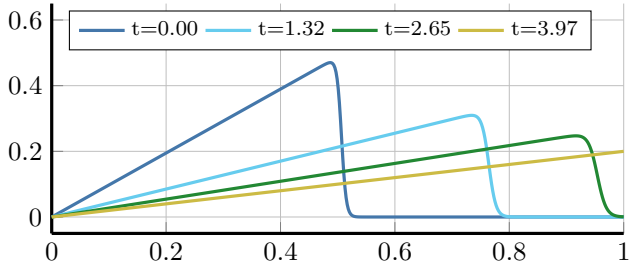

As the first example, we consider the one-dimensional viscous Burgers' equation
\begin{align*}
	\partial_t \state(t,\xi) = \tfrac{1}{\mathrm{Re}} \partial_{\xi^2} \state(t,\xi) - \state(t,\xi)\partial_\xi \state(t,\xi),
\end{align*}
with known analytical solution
\begin{equation}
	\label{eqn:BurgersSol}
	\state(t,\xi) = \frac{\xi}{t+1}\left(1+\sqrt{\frac{t+1}{\exp\left(\tfrac{\mathrm{Re}}{8}\right)}}\exp\left(\mathrm{Re}\tfrac{\xi^2}{4t+4}\right)\right)^{-1},
\end{equation}
taken from \cite{MauLB21} (see also \cite{BlaSU22}, where a similar setting is used).
The spatial domain is set to $\Omega \vcentcolon= [0,1]$, the time interval to $\timeInt \vcentcolon= [0,4]$, and the Reynolds number is given by $\mathrm{Re} = 1000$. 
For an illustration of the solution of the full-order model, we refer to Figure~\ref{fig:BurgersData}. 

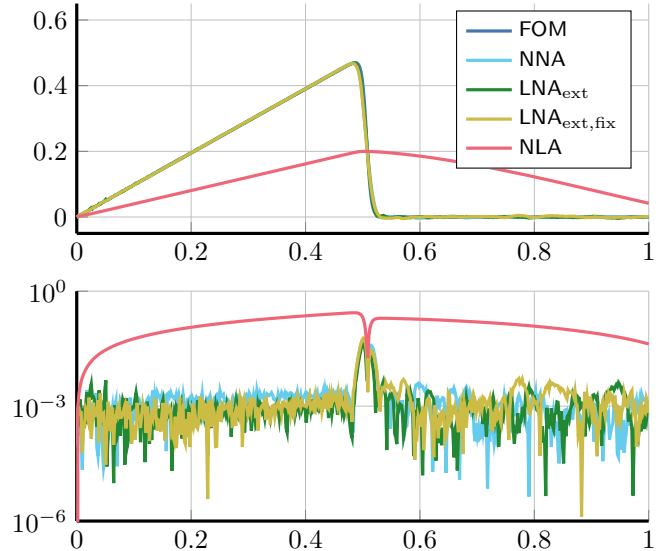
\begin{figure}
	\begin{flushright}
   \begin{tikzpicture}
    \begin{axis}[
    width=.85\linewidth,
    height=1.2in,
    scale only axis,
    xmin=0, xmax=1,
    ymin=-0.05, ymax=0.65,
	grid=both,
	grid style={line width=.1pt, draw=gray!10},
	major grid style={line width=.2pt,draw=gray!50},
	axis lines*=left,
	axis line style={line width=\lineWidth},
    legend entries={\FOM, \NNA, \LNAwt, \LNAextTime, \NLA},
    legend style={font=\footnotesize},
    legend cell align={left},
    legend pos = north east,
    cycle list name=approx,
    no markers,
    ]
    \addplot table[x=x, y=FOM,col sep=comma] {BurgersEquationLargeTime-approx-t-0.03.dat};
    \addplot table[x=x, y=nonlin-N512-r1-f3,col sep=comma] {BurgersEquationLargeTime-approx-t-0.03.dat};
    \addplot table[x=x, y=lin_nonlin-N513-r1-f3,col sep=comma] {BurgersEquationLargeTime-approx-t-0.03.dat};
    \addplot table[x=x, y=linOnlyTime_nonlin-N513-r1-f3,col sep=comma] {BurgersEquationLargeTime-approx-t-0.03.dat};
    \addplot table[x=x, y=nonlin_lin-N512-r1-f3,col sep=comma] {BurgersEquationLargeTime-approx-t-0.03.dat};
\end{axis}
\end{tikzpicture}

\begin{tikzpicture}
    \begin{axis}[
    width=.85\linewidth,
    height=1.2in,
    scale only axis,
    xmin=0, xmax=1,
    ymin=1e-6, ymax=1,
    ymode={log},
	grid=both,
	grid style={line width=.1pt, draw=gray!10},
	major grid style={line width=.2pt,draw=gray!50},
	axis lines*=left,
	axis line style={line width=\lineWidth},
    cycle list name=error,
    no markers,
    ]
    \addplot table[x=x, y=nonlin-N512-r1-f3,col sep=comma] {BurgersEquationLargeTime-error-t-0.03.dat};
    \addplot table[x=x, y=lin_nonlin-N513-r1-f3,col sep=comma] {BurgersEquationLargeTime-error-t-0.03.dat};
    \addplot table[x=x, y=linOnlyTime_nonlin-N513-r1-f3,col sep=comma] {BurgersEquationLargeTime-error-t-0.03.dat};
    \addplot table[x=x, y=nonlin_lin-N512-r1-f3,col sep=comma] {BurgersEquationLargeTime-error-t-0.03.dat};
\end{axis}
\end{tikzpicture}
	\end{flushright}
	\caption{Burgers' equation: approximations (top) and corresponding errors in absolute value (bottom) for the solution at $t=0.03$ over different autoencoders with $\stateDimRed=1$ and scenario (A). The first 4 lines in the top figure nearly coincide.}
	\label{fig:Burgers:approx}
\end{figure}

\begin{figure}
	\begin{flushright}
	\begin{tikzpicture}
    \begin{axis}[
    width=.85\linewidth,
    height=1.2in,
    scale only axis,
    xmin=0, xmax=1,
    ymin=-0.05, ymax=0.4,
	grid=both,
	grid style={line width=.1pt, draw=gray!10},
	major grid style={line width=.2pt,draw=gray!50},
	axis lines*=left,
	axis line style={line width=\lineWidth},
    legend columns=4,
    legend entries={\POD $t_1$, \NLA $t_1$,\POD $t_2$, \NLA $t_2$},
    legend style={font=\footnotesize},
    legend cell align={left},
    legend pos = north west,
    cycle list name=approxEverySecondDashed
    ]
    \addplot table[x=sampling, y=POD-t-1.35,col sep=comma] {BurgersEquationLargeTime-POD-approx.dat};
    \addplot table[x=x, y=nonlin_lin-N512-r1-f3,col sep=comma] {BurgersEquationLargeTime-approx-t-1.35.dat};    
    \addplot table[x=sampling, y=POD-t-4.00,col sep=comma] {BurgersEquationLargeTime-POD-approx.dat};
    \addplot table[x=x, y=nonlin_lin-N512-r1-f3,col sep=comma] {BurgersEquationLargeTime-approx-t-4.00.dat};    
\end{axis}
\end{tikzpicture}
	\end{flushright}
	\caption{Burgers' equation: \POD (solid) and \NLA (dashed) approximation for times $t_1= 1.35$, $t_2=4.00$.} 
	\label{fig:Burgers:PODNLA}
\end{figure}
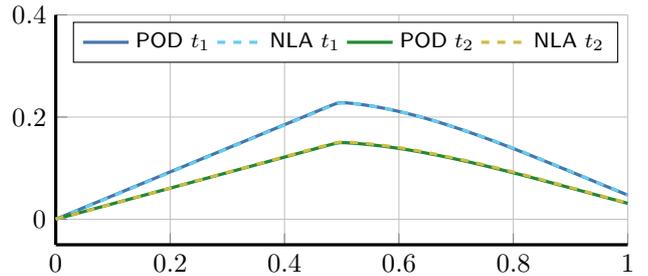

The error on the testing data for the best run is reported in Table~\ref{tab:Burgers:bestTestingErr}. The highlighted cell corresponds to the minimal value in the column.
\begin{table}[ht]
	\centering
	\caption{Burgers' equation: Best relative testing error with respect to \eqref{eqn:error_measure}.}
	\label{tab:Burgers:bestTestingErr}
	\begin{tabular}{lrrr}
		\toprule 
		method & (A) & (B) & (C) \\\midrule
		\NNA & \num{0.011006} & \num{0.011094} & \num{0.013203}\\
		\LNAext &  \cellcolor{cellgray}{\num{0.009064}} &  \cellcolor{cellgray}{\num{0.010766}} &  \cellcolor{cellgray}{\num{0.012751}}\\
		\LNA &{\num{0.011670}} & \num{0.013307}  & \num{0.020180}\\
		\LNAextTime & \num{0.009687} & {\num{0.011000}}  & {\num{0.015159}}\\
		\NLA & \num{0.491939} & \num{0.491762}  & \num{0.491762}\\\bottomrule
	\end{tabular}
\end{table}
We report that the best relative testing errors are achieved by the \LNAext, with the second-best error obtained by either the standard \NNA or the \LNAextTime. Further, we observe that the time-extended data set is preferable in this example for linear-nonlinear autoencoder configurations to obtain a good approximation error, as the \LNA obtains worse results than the standard \NNA. As expected from Theorem~\ref{thm:linearDecoder}, the configuration approximation error of the \NLA is bounded below by the Kolmogorov barrier. Thus, it is not expected to yield better results than a linear-linear configuration. Indeed, the \NLA yields relative errors that almost precisely match the relative error of \POD with reduced dimension $\stateDimRed =1$, i.e., \POD-error \num{0.4927}. To achieve a relative error of less than 2\% with \POD, we must use a reduced dimension of at least $\stateDimRed = 21$. For visualization, we plot the different approximations and their errors for the time $t=1.35$ in Figure~\ref{fig:Burgers:approx}. To further illustrate how similar the approximation results of the \NLA and the \POD are, we plot both in Figure~\ref{fig:Burgers:PODNLA} for two separate time instances.

To further inspect the training of the different methods, we also report the average testing error of all 100 model runs in Table~\ref{tab:Burgers:averageTestingErr}. 
\begin{table}[ht]
	\centering
	\caption{Burgers' equation: average relative testing error with respect to \eqref{eqn:error_measure}.}
	\label{tab:Burgers:averageTestingErr}
	\begin{tabular}{lrrr}
		\toprule 
		method & (A) & (B) & (C) \\\midrule
		\NNA & \num{0.037293} & \num{0.024301} & \num{0.026543}\\
		\LNAext & {\num{0.016175}} & {\num{0.015876}} &  \cellcolor{cellgray}{\num{0.020466}}\\
		\LNA  &{\num{0.024102}} & \num{0.024756}  & \num{0.028304}\\
		\LNAextTime &  \cellcolor{cellgray}{\num{0.016136}} &  \cellcolor{cellgray}{\num{0.015617}}  & {\num{0.021201}}\\
		\NLA & \num{0.492049} & \num{0.492074}  & \num{0.492155}\\\bottomrule
	\end{tabular}
\end{table}
We observe that the best average testing error is either obtained by the \LNAextTime$\,$or by the \LNAext, which thus yield the most stable configurations. So again, the lowest errors are reported by configurations utilizing the additional time information in training, where the simpler configuration \LNAextTime$\,$outperforms the \LNAext$\,$in two instances. On the other hand, autoencoders trained on the normal data set increase by a higher percentage, rendering the configurations less robust. Again, as expected, the \NLA yields robust results, as it mimicks a linear-linear model.

\begin{table}[ht]
	\centering
	\caption{Burgers' equation: average training time per epoch.}
	\label{tab:Burgers:averageTrainingTime}
	\begin{tabular}{lrrr}
		\toprule 
		method & (A) & (B) & (C) \\\midrule
		\NNA & \num{0.040985} & \num{0.034736} & \num{0.029256}\\	
		\LNAext & {\num{0.022695}} & {\num{0.018632}} & {\num{0.016089}}\\
		\LNA &{\num{0.021361}} & \num{0.018025}  & \num{0.017738}\\			
		\LNAextTime&  \cellcolor{cellgray}{ \num{0.020319}} &  \cellcolor{cellgray}{\num{0.017131}}  &  \cellcolor{cellgray}{\num{0.015157}}\\ \bottomrule
	\end{tabular}
\end{table}

In Table~\ref{tab:Burgers:averageTrainingTime}, we compare the average training time per epoch and consider all configurations with a nonlinear decoder. First, note that using a linear encoder pays off in the average training time, as all linear-nonlinear autoencoder configurations attain a training time of 60\% or less compared to the \NNA. Further, we observe that the \LNAextTime$\,$achieves the fastest times, which we assume is due to the even simpler design in the encoder. 
Thus, we achieve a similar or better approximation error with a linear-nonlinear autoencoder with time-extended data while needing less average training time per epoch. 

\subsection{Advection equation}
\label{subsec:advection}

For our second example, we use the advection equation
\begin{align*}
	\partial_t  \state(t,\xi) + c\partial_\xi \state(t,\xi) = 0,
\end{align*}
with sawtooth initial value
\begin{align*}
	\state_0(\xi) =  \begin{cases} 
		\tfrac{1}{\sigma} (\xi - \beta), & \beta \le \xi \le \beta+\sigma, \\
		0, & \textnormal{else},
	\end{cases}
\end{align*}
such that the solution is given by $\state(t,\xi) = \state(\xi - ct)$. For the computational setup, we have $\Omega = [0,1]$, $\timeInt = [0,1]$, $c=1$, $\sigma = 0.1$ and $\beta = 0$. Note that with this choice, the time-to-solution map is not smooth, and hence, it does not satisfy the assumptions from Theorems~\ref{thm:cor3-6} and~\ref{thm:linEncoder}. Nevertheless, despite the assumptions not being satisfied, we achieve similar results as in the previous example. 

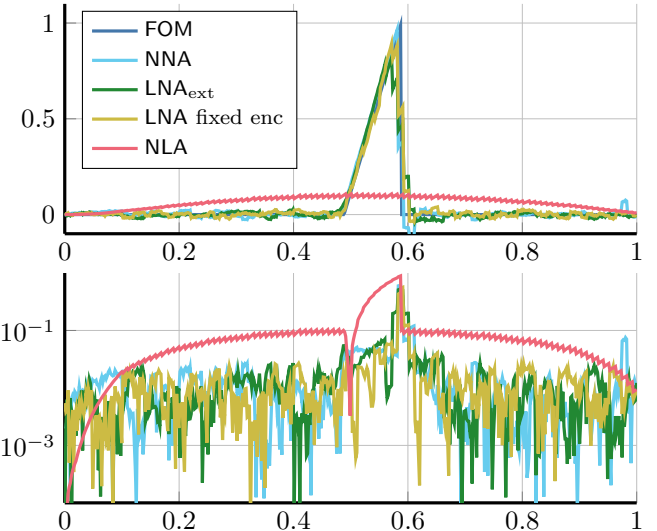
\begin{figure}
	\begin{flushright}
	\begin{tikzpicture}
    \begin{axis}[
    width=.85\linewidth,
    height=1.2in,
    scale only axis,
    xmin=0, xmax=1,
    ymin=-0.1, ymax=1.1,
	grid=both,
	grid style={line width=.1pt, draw=gray!10},
	major grid style={line width=.2pt,draw=gray!50},
	axis lines*=left,
	axis line style={line width=\lineWidth},
    legend entries={\FOM, \NNA, \LNAwt, \LNA fixed enc, \NLA},
    legend style={font=\footnotesize},
    legend cell align={left},
    legend pos = north west,
    cycle list name=approx,
    no markers,
    ]
    \addplot table[x=x, y=FOM,col sep=comma] {sawTooth-approx-t-0.49.dat};
    \addplot table[x=x, y=nonlin-N512-r1-f3,col sep=comma] {sawTooth-approx-t-0.49.dat};
    \addplot table[x=x, y=lin_nonlin-N513-r1-f3,col sep=comma] {sawTooth-approx-t-0.49.dat};
    \addplot table[x=x, y=linOnlyTime_nonlin-N513-r1-f3,col sep=comma] {sawTooth-approx-t-0.49.dat};
    \addplot table[x=x, y=nonlin_lin-N512-r1-f3,col sep=comma] {sawTooth-approx-t-0.49.dat};
\end{axis}
\end{tikzpicture}

\begin{tikzpicture}
    \begin{axis}[
    width=.85\linewidth,
    height=1.2in,
    scale only axis,
    xmin=0, xmax=1,
    ymin=1e-4, ymax=1,
    ymode={log},
	grid=both,
	grid style={line width=.1pt, draw=gray!10},
	major grid style={line width=.2pt,draw=gray!50},
	axis lines*=left,
	axis line style={line width=\lineWidth},
    cycle list name=error,
    no markers,
    ]
    \addplot table[x=x, y=nonlin-N512-r1-f3,col sep=comma] {sawTooth-error-t-0.49.dat};
    \addplot table[x=x, y=lin_nonlin-N513-r1-f3,col sep=comma] {sawTooth-error-t-0.49.dat};
    \addplot table[x=x, y=linOnlyTime_nonlin-N513-r1-f3,col sep=comma] {sawTooth-error-t-0.49.dat};
    \addplot table[x=x, y=nonlin_lin-N512-r1-f3,col sep=comma] {sawTooth-error-t-0.49.dat};
\end{axis}
\end{tikzpicture}
	\end{flushright}
	\caption{Advection equation: approximations (top) and corresponding errors in absolute values (bottom) for the solution at $t=0.49$ over different autoencoders with $\stateDimRed=1$ and scenario (A).}
	\label{fig:advectionEquation:approx}
\end{figure}

First, we report the best error in Table~\ref{tab:Advection:bestTestingErr}. Although not covered by our theory, we observe that linear-nonlinear autoencoder configurations achieve the best errors in agreement with our first experiment. In particular, \LNAextTime$\,$performs best over all parametrizations, with \NNA yielding the second-best results. In contrast, \NLA is, again as expected, significantly worse; its error matches the relative error of \POD with reduced dimension $\stateDimRed =1$ (\POD-error \num{0.9615}). For a relative error of less than 30\% with \POD, we need a reduced dimension of $\stateDimRed = 35$ in this example. For visualization we plot the different approximations and their errors for the time $t=0.49$ in Figure~\ref{fig:advectionEquation:approx}.
\begin{table}[ht]
	\centering
	\caption{Advection equation: best relative testing error with respect to \eqref{eqn:error_measure}.}
	\label{tab:Advection:bestTestingErr}
	\begin{tabular}{lrrrr}
		\toprule 
		method & (A) & (B) & (C)\\\midrule
		\NNA & \num{0.366163} & \num{0.391770} & \num{0.379639}\\
		\LNAext & \num{0.380621} & \num{0.404984} & \num{0.400513}\\
		\LNA  & \num{0.392716} & \num{0.390943} & \num{0.407499}\\
		\LNAextTime & \cellcolor{cellgray}{\num{0.282484}} & \cellcolor{cellgray}{\num{0.296808}} & \cellcolor{cellgray}{\num{0.323463}}\\
		\NLA & \num{0.961535} & \num{0.961471} & \num{0.961633}\\\bottomrule
	\end{tabular}
\end{table}

The best average testing errors are reported in Table~\ref{tab:Advection:AvgTestingErr}. The \LNAextTime$\,$yields the most robust configuration, as observed in the previous example. We notice that the percentage increase is lowest for the \LNAextTime$\,$with 2\%-16\% compared to the other configurations, with 29\%-50\% increase. We excluded the \NLA because, despite its stability, it has a significantly higher error than other configurations due to the Kolmogorov barrier.

\begin{table}[ht]
	\centering
	\caption{Advection equation: average relative testing error with respect to \eqref{eqn:error_measure}.}
	\label{tab:Advection:AvgTestingErr}
	\begin{tabular}{lrrrr}
		\toprule 
		method & (A) & (B) & (C)\\\midrule
		\NNA & \num{0.551107} & \num{0.503646} & \num{0.514109}\\
		\LNAext & \num{0.567999} & \num{0.586701} & \num{0.599578}\\
		\LNA  & \num{0.564725} & \num{0.565638} & \num{0.607045}\\
		\LNAextTime & \cellcolor{cellgray}{\num{0.328802}} & \cellcolor{cellgray}{\num{0.333614}} & \cellcolor{cellgray}{\num{0.372758}}\\
		\NLA & \num{0.961540} & \num{0.961503} & \num{0.961618}\\\bottomrule
	\end{tabular}
\end{table}

\begin{figure}
	\begin{flushright}
	\begin{tikzpicture}
    \begin{axis}[
    width=.85\linewidth,
    height=1.2in,
    scale only axis,
    xmin=0, xmax=1,
    ymin=1e-4, ymax=1.1,
    ymode={log},
	grid=both,
	grid style={line width=.1pt, draw=gray!10},
	major grid style={line width=.2pt,draw=gray!50},
	axis lines*=left,
	axis line style={line width=\lineWidth},
    legend columns=4,
    legend entries={\NNA, \LNAext, \LNA, \LNAextTime},
    legend style={font=\footnotesize},
    legend cell align={left},
    legend pos = north west,
    cycle list name=error,
    no markers,
    ]
    \addplot table[x=sampling, y=projectionProp,col sep=comma] {sawTooth-projProp-nonlin-N512-r1-f3.dat};
    \addplot table[x=sampling, y=projectionProp,col sep=comma] {sawTooth-projProp-lin_nonlin-N513-r1-f3.dat};
     \addplot table[x=sampling, y=projectionProp,col sep=comma] {sawTooth-projProp-lin_nonlin-N512-r1-f3.dat};
    \addplot table[x=sampling, y=projectionProp,col sep=comma] {sawTooth-projProp-linOnlyTime_nonlin-N513-r1-f3.dat};
\end{axis}
\end{tikzpicture}
	\end{flushright}
	\caption{Advection equation: point projection property~\eqref{eqn:projectionProperty} for the trained autoencoders \NNA, \LNAext, \LNA, \LNAextTime$\,$with $\stateDimRed=1$ and scenario~(A).}
	\label{fig:advectionEquation:projectionProperty}
\end{figure}
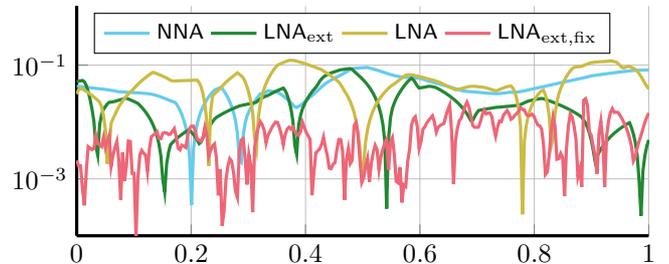

Finally, we investigate the point projection property~\eqref{eqn:projectionProperty} for the different autoencoders. We emphasize that neither of the autoencoders strictly enforces the point projection property during training, nor a deviation of it is penalized. However, in \cite[Thm.~6.4]{BucGHU24}, it is shown that the error of the point projection can be bounded by the least squares error in training. In particular, we expect a smaller deviation from the point projection property for configurations with smaller validation errors. We investigate this numerically for the autoencoder configurations \NNA, \LNAext, \LNA, \LNAextTime$\,$in Figure~\ref{fig:advectionEquation:projectionProperty}, where we plot the deviation from the point projection property (in absolute value). 
We observe that even if this property was not enforced directly in the training, the \LNAextTime$\,$obtains the smallest deviation, which is reasonable since it had the best average approximation error. 

\section{Conclusion}\label{sec:conclusion}

We present a novel approach for autoencoder training for \MOR by incorporating time and parameter information into the training data. We show that a linear encoder suffices to achieve the same accuracy as the \NNA, illustrated in two proof-of-concept examples. For the Burgers' equation, the \LNAext$\,$matches the \NNA's accuracy while having faster training times per epoch. The second example with the advection equation shows that similar results can be obtained even when theoretical assumptions are violated. Additionally, we assess the quality of the point projection property, which is not explicitly included in training. Nevertheless, the \LNAextTime$\,$configuration yields stable results consistent with the low training error.

With these proof-of-concept examples, the next step is to numerically validate the theory on further benchmarks and for different autoencoder configurations, such as convolutional autoencoders. Further, we aim to extend this work to the parametric setting and examine how the \ROM error relates to the projection error. Eventually, we plan to use the resulting outcomes for the efficient (optimal) control of transport-dominated partial different equations, e.g., extending the recent work by \cite{BrBS2024}. 

\bibliography{journalAbbr.bib,literature.bib}
                                                   
\end{document}